\documentclass[a4paper,11pt]{amsart}

\usepackage{amsthm, amsmath, latexsym, amsfonts, epsfig}
\usepackage[all]{xy}
\usepackage{amssymb}

\newtheorem{Lemma}{Lemma}[section]
\newtheorem{Proposition}[Lemma]{Proposition}
\newtheorem{Theorem}[Lemma]{Theorem}
\newtheorem{remark}[Lemma]{{\em Remark}}
\newtheorem{Corollary}[Lemma]{Corollary}

\newcommand{\Pic}{\mathrm{Pic}}

\newcommand{\Aut}{\mathrm{Aut}}

\newcommand{\Sp}{\mathrm{Sp}}
\newcommand{\Z}{\mathbb{Z}}
\newcommand{\Q}{\mathbb{Q}}

\renewcommand{\P}{\mathbb{P}}
\newcommand{\mg}{\mathcal{M}_g}

\newcommand{\mtwobar}{\overline{\mathcal{M}}_2}

\newcommand{\mobar}{\overline{\mathcal{M}}_{0,n}}
\newcommand{\sgbar}{\overline{S}_g}
\newcommand{\sgnbar}{\overline{S}_{g,n}}
\newcommand{\stwobar}{\overline{S}_2}
\newcommand{\sgnm}{\overline{S}_{g,n}^{\hspace{0.05cm}(m_1, \ldots, m_n)}}
\newcommand{\spiu}{\overline{S}_2^+}
\newcommand{\smeno}{\overline{S}_2^-}

\title{On the geometry of $\stwobar$} \author{Gilberto Bini and
Claudio Fontanari}

\email{gilberto.bini@mat.unimi.it} \curraddr{{\sc Dipartimento di
Matematica \\ Universit\`a degli Studi di Milano \\ Via C. Saldini 50 \\
20133 Milano \\ Italy.}}

\email{claudio.fontanari@polito.it}\curraddr{
{\sc Dipartimento di Matematica \\ Politecnico di Torino \\
Corso Duca degli Abruzzi 24 \\ 10129 Torino \\ Italy.}}


\begin{document}

\begin{abstract}
We investigate topological properties of the moduli space of spin structures over genus two curves. In particular, we provide a combinatorial description of this space and give a presentation of the (rational) cohomology ring via generators and relations. 

\emph{The authors warn the reader that an erratum to the present paper has been kindly posted on arXiv by Sebastian Krug.}

\end{abstract}

\maketitle

\section{Introduction}\label{sec1}

The moduli space $\sgbar$ of spin curves of genus $g$ has been introduced in \cite{Cornalba:89} in order to compactify the moduli space of pairs 
$$
\left(\hbox{smooth genus}\; g \; \hbox{complex curve}\; C, \; \hbox{theta-characteristic on}\; C \right).
$$ 

The interest in this space has been increasing in the last few years
due to a wide spectrum of applications, which range from Mathematical
Physics to Arithmetic Geometry passing through Algebraic Geometry and
Combinatorics (see, for instance, \cite{AbrJar:01}, \cite{BF},
\cite{CaCa}, \cite{CaCCo}, \cite{Chi}, \cite{Chie}, \cite{JKV},
\cite{Pacs}).

Nonetheless, very little is known about the topological properties of $\sgbar$. The first results in this direction can be found in \cite{Harer:90} and \cite{Harer:93}. In these pioneering articles, the stability of the homology groups and some low degree cohomology groups of $S_g$ are investigated.  

In the present paper, we compute the virtual cohomological dimension of $S_g$ (see Lemma ~\ref{virtual}). Next, we focus on the compactification $\sgbar$, especially on the first relevant case $g=2$. In particular, we determine the Hodge diamond of ${\overline S}_2$ (see Proposition \ref{hodge}) and describe the cohomology ring of  ${\overline S}_2$ via generators and relations (see Theorem \ref{coho2n} and Theorem \ref{coho3n}).

At first sight, it might seem that the topology of the moduli space of
spin curves is not really different from that of the moduli space of
stable curves. In fact, the constructions of these two spaces do not
differ in a substantial way. On the other hand, $\sgbar$ turns out to
have a richer geometric structure, which requires a more careful
analysis. Indeed, $\sgbar$ is the union of two connected components,
$\sgbar^+$ and $\sgbar^-$, which correspond to even and odd spin
curves, respectively. In the case $g=2$, we separately investigate the
cohomology of each component and  discover two non-isomorphic ring
structures.

In order to accomplish this task, we combine an inductive approach
inspired by \cite{ArbCor:98} with direct calculations in the style of
\cite{Mum}. This leads us to introduce suitable moduli spaces of
pointed spin curves and to address several related problems. In 
particular, we obtain an explicit affine stratification of
${\overline S}_2$ (see Appendix A).

Thanks are due to Maurizio Cornalba for helpful conversations on the
subject of this research and to Katharina Ludwig for pointing out a
mistake in a previous version of the present paper. The authors are
deeply grateful to the referee for careful reading of the manuscript 
and patient advice on how to improve it 
(see for instance Remark~\ref{referee}). 

\emph{An erratum to the present paper has been kindly posted on arXiv by Sebastian 
Krug (see \cite{krug}).}

Throughout we work over the field ${\mathbb C}$ of complex numbers.

\section{The moduli space of curves with spin structures}\label{sec2}

In this section, we recall some basic definitions about spin structures. Here we follow closely \cite{Cornalba:89}, to which we refer for more details. For a more general approach see, for instance, \cite{AbrJar:01}.

Let $X$ be a Deligne-Mumford semistable curve and let $E$ be a
complete, irreducible subcurve of $X$. The curve $E$ is said to be
\emph{exceptional} when it is smooth, rational, and intersects the other
components in exactly two points. Moreover, $X$ is said to be
\emph{quasi-stable} when any two distinct exceptional components of $C$
are disjoint. In the sequel, $\tilde{X}$ will denote the subcurve
$\overline{X \setminus \cup E_i}$ obtained from $X$ by removing all the
exceptional components.

A \emph{spin curve} of genus $g$ (see \cite{Cornalba:89}, \S~2)
is the datum of a quasi-stable genus $g$ curve $X$ with an invertible 
sheaf $\zeta_X$ of degree $g-1$ on $X$ and a homomorphism of invertible 
sheaves $
\alpha_X: \zeta_X^{\otimes 2} \longrightarrow \omega_X
$
such that  i) $\zeta_X$ has degree $1$ on every exceptional component of $X$, and ii) $\alpha_X$ is not zero at a general point of every non-exceptional 
component of $X$. Therefore, $\alpha_X$ vanishes identically 
on all exceptional components of $X$ and induces an isomorphism 
$
\tilde{\alpha}_X: \zeta_X^{\otimes 2} \vert_{\tilde{X}} \longrightarrow 
\omega_{\tilde{X}}.
$
In particular, when $X$ is smooth, $\zeta_X$ is just a theta-characteristic 
on $X$.  Two spin curves 
$(X, \zeta_X, \alpha_X)$ and $(X', \zeta_{X'}, \alpha_{X'})$ are 
\emph{isomorphic} if there are isomorphisms $\sigma: X \to X'$ and 
$\tau: \sigma^*(\zeta_{X'}) \to \zeta_X$
such that $\tau$ is compatible with the natural 
isomorphism between $\sigma^*(\omega_{X'})$ and $\omega_X$.

A \emph{family of spin curves} is a flat family of quasi-stable 
curves $f: \mathcal{X} \to S$ with an invertible sheaf $\zeta_f$ on 
$\mathcal{X}$ and a homomorphism
$
\alpha_f: \zeta_f^{\otimes 2} \longrightarrow \omega_f
$
such that the restriction of these data to any fiber of $f$ gives rise 
to a spin curve.

\noindent Two families of spin curves $f: \mathcal{X} \to S$ and 
$f': \mathcal{X}' \to S$ are \emph{isomorphic} if there are isomorphisms
$\sigma: \mathcal{X} \longrightarrow \mathcal{X}'$ and 
$\tau: \sigma^*(\zeta_{f'}) \longrightarrow \zeta_f$
such that $f = f' \circ \sigma$ and $\tau$ is compatible with the natural 
isomorphism between $\sigma^*(\omega_{f'})$ and $\omega_f$. 

\noindent Let $\sgbar$ be the moduli space of isomorphism classes of spin curves of 
genus $g$. Denote by $S_g$ the open subset consisting of classes of smooth curves. 
As shown in \cite{Cornalba:89}, \S \; 5, $\sgbar$ has  a natural structure of analytic orbifold given as follows. For any spin curve $X$,  
there is a neighbourhood $U$ of $[X]$ such that
$U \cong B_X / \Aut(X)$, where $B_X$ is a $3g-3$-dimensional polydisk and $\Aut(X)$ is  the automorphism group of the spin curve $X$. Alternatively, ${\overline S}_g$ may be viewed as a projective normal variety with finite quotient singulaties.

The moduli space of spin curves can be slightly generalized 
as follows. For all integers $g$, $n$, $m_1, \ldots, m_n$, 
such that $2g-2+n>0$, $0 \le m_i \le 1$ for every $i$, and 
$\sum_{i=1}^n m_i$ is even, we define

\begin{eqnarray*}
\sgnm &:=& \Big \{ [(C, p_1, \ldots, p_n; \zeta; \alpha)]: 
(C, p_1, \ldots, p_n) \textrm{ is a genus $g$}  \\ 
& & \textrm{quasi-stable projective curve with $n$ marked points}; \\
& &\zeta \textrm{ is a line bundle of degree $g-1 + \frac12\sum_{i=1}^n m_i$ on $C$ having} \\
& & \textrm{ degree $1$ on every exceptional component of $C$, and} \\ 
& &\alpha: \zeta^{\otimes 2} \to \omega_C(\sum_{i=1}^n m_i p_i) 
\textrm{ is a homomorphism which} \\ 
& & \textrm{is not zero at a general point of every 
non-exceptional}\\ 
& & \textrm{component of $C$} \Big \}.
\end{eqnarray*}

In order to put an analytic structure on $\sgnm$, we notice that
Cornalba's construction in \cite{Cornalba:89} can be easily adapted to
$\sgnm$. Indeed, from the universal deformation of the stable model of
$(C, p_1, \ldots, p_n)$ we obtain exactly as in \cite{Cornalba:89},
\S~4, a universal deformation $\mathcal{U}_X \to B_X$ of $X = (C, p_1,
\ldots, p_n; \zeta; \alpha)$. Next, we put on $\sgnm$ the structure of of the quotient analytic space $B_X /$ $\Aut(X)$ following \cite{Cornalba:89}, \S~5. Alternatively,
we can regard $\sgnm$ as the coarse moduli space associated to the
stack of $r$-spin curves (in the easiest case $r=2$), which has been
constructed by Jarvis in \cite{Jarvis:00} and revisited by Abramovich
and Jarvis in \cite{AbrJar:01}.

Analogously to $\sgbar$ (see \cite{Cornalba:89}, Proposition~5.2),
the spaces $\sgnm$ are normal projective varieties of complex dimension 
$3g-3+n$. If $m_1 = \ldots = m_n = 0$, then $\sgnbar := 
\sgnbar^{\hspace{0.05cm}(0, \ldots, 0)}$ splits into two disjoint 
irreducible components $\sgnbar^+$ and $\sgnbar^-$ that consist of 
the even and the odd spin curves, respectively (see \cite{Cornalba:89}, Lemma~6.3). Similarly to Lemma 1 in \cite{BF}, we point out the following

\begin{Proposition}\label{iso}
Let $\Pic(\sgnm) := H^1(\sgnm, \mathcal{O}^*)$.
There is a natural isomorphism 
$$ 
\Pic(\sgnm) \otimes \Q \stackrel{\cong}{\longrightarrow}
A_{3g-4+n}(\sgnm) \otimes \Q.
$$
\end{Proposition}
 
\proof Since $\sgnm$ is normal, there is an injection:
$$
\Pic(\sgnm) \otimes {\mathbb Q} \hookrightarrow A_{3g-4+n}(\sgnm) \otimes {\mathbb Q}.
$$
Moreover, from the construction of $\sgnm$ it follows that the singularities of $\sgnm$ are of finite quotient type, so every Weil divisor is ${\mathbb Q}$-Cartier and there is a surjective morphism:
$$
\Pic(\sgnm) \otimes {\mathbb Q} \rightarrow A_{3g-4+n}(\sgnm) \otimes {\mathbb Q}.
$$

Hence the claim follows.
\qed

We end this section with some results on the topology of $S_g$, which will be applied in the next Section when $g=2$. In what follows, all homology and cohomology groups are intended to have rational coefficients.

Let us recall some general notions from \cite{Harer:86}, \S~4.
A group $\Gamma$ is \emph{virtually torsion-free} if it has a subgroup
$G$ of finite index which is torsion free. If $\Gamma$ is virtually
torsion free, the \emph{virtual cohomological dimension} of $\Gamma$
is the cohomological dimension of a torsion free subgroup $G$ of
finite index in $\Gamma$. A theorem of Serre states that this number
is independent of the choice of $G$ (see \cite{Harer:86}, p.~173).

It is well-known that the mapping class group $\Gamma_g$ is virtually 
torsion-free (see for instance \cite{Harer:86}, p.~172). Moreover, 
Harer proved that the virtual cohomological dimension of 
$\Gamma_g$ is $4g-5$ (see \cite{Harer:86}, Theorem~4.1).
Since the Teichm\"uller space $T_g$ is contractible, the rational 
homology of $\mg$ and $\Gamma_g$ is the same, so the group $H_k(\mg, \Q)$ 
is zero for $k > 4g - 5$ (see \cite{Harer:86}, Corollary~4.3).

Let now $S$ be a smooth orientable surface of genus $g$ and let $Q$ 
be a $\Z / 2 \Z$-quadratic form on $H_1(S,\Z / 2 \Z)$. The isomorphism
class of $Q$ is determined by its Arf invariant $\varepsilon = 0$ or $1$
(see \cite{Harer:90}, p. 324). Call $Q_0$ (resp. $Q_1$)
the even (resp. the odd) quadratic form and let $G_g(Q_i)$ be the 
subgroup of the mapping class group $\Gamma_g$ which preserves $Q_i$. 
Then the following holds:
 
\begin{Lemma}\label{virtual} 
For $i = 0, 1$, $G_g(Q_i)$ is virtually torsion-free and its virtual 
cohomological dimension is $4g-5$.
\end{Lemma}   

\proof Let $\overline{\chi}: \Gamma_g \longrightarrow \Sp(2g, \Z / 2
\Z)$ be the following homomorphism. The image under $\overline{\chi}$
of an isotopy class $[\gamma]$ of orientation-preserving diffeomorphism of $S$
is the automorphism induced by  $\gamma$ on $H_1(S, \Z / 2 \Z)$ (see \cite{ACGH:02},
Chapter~14, (2.1)).  Since the target group is finite, the kernel $K$
of $\overline{\chi}$ has finite index in $\Gamma_g$. It is a
straightforward consequence of the definitions that $K \subseteq
G_g(Q_i)$, hence we may deduce that $G_g(Q_i)$ has finite index in
$\Gamma_g$. Now, pick a torsion-free subgroup $H$ of $\Gamma_g$ of
finite index and consider $G_i := H \cap G_g(Q_i)$.  Since both $H$
and $G_g(Q_i)$ have finite index in $\Gamma_g$, it follows that $G_i$
has finite index in $\Gamma_g$. Since $H$ is torsion-free, we have
that $G_i$ is torsion-free. As mentioned before, all torsion-free
subgroups of finite index in $\Gamma_g$ have the same cohomological
dimension. Therefore,  from the previous two facts and Harer's theorem
we may deduce that the cohomological dimension of $G_i$ is $4g-5$.
Since $G_i$ has finite index in $\Gamma_g$, it has finite
index in $G_g(Q_i)$ too. Thus, it computes the virtual cohomological
dimension of $G_g(Q_i)$.

\qed
 
Recall that $S_g$ is the disjoint union of $S_g^+$, which contains 
curves with even theta-characteristics, and $S_g^-$, which contains those with 
odd theta-characteristics. Further, we have $S_g^+ = T_g / G_g(Q_0)$ and 
$S_g^- = T_g / G_g(Q_1)$ (see \cite{Harer:90}, p. 324). Analogously to $\mg$, 
the following holds.

\begin{Theorem}\label{harer} We have
$H_k(S_g) = 0$ for $k > 4g-5$.
\end{Theorem}

\section{The Rational cohomology of $\stwobar$}\label{sec6}

In this section, we specialize to the case $g=2$ and we investigate the 
additive and multiplicative structure of the cohomology algebra of 
${\overline S}_2$. Let $\nu: \stwobar \rightarrow \mtwobar$ be the 
morphism which forgets the spin structure and passes
to the stable model. It is a $16:1$ ramified covering of normal
threefolds. Denote by $\nu^+: \spiu \rightarrow \mtwobar$
(resp. $\nu^-: \smeno \rightarrow \mtwobar$) the restriction of $\nu$
to the moduli space of curves with even spin structures $\spiu$
(resp. odd spin structures $\smeno$).

\begin{remark}\label{remaff}
The stratification of ${\overline S}_2$ by topological type is affine.
Indeed, as shown in \cite{Mum}, the stratification of $\mtwobar$ by 
graph type is affine and since $\nu$ is finite and preserving topological 
type, ${\overline S}_2$ turns out to be stratified as the union of affine 
subvarieties.
\end{remark}

For an explicit description of the stratification in Remark
\ref{remaff}, see Appendix A.

\begin{Theorem}\label{odd}
We have $H^1(\stwobar)=H^3(\stwobar)=0$.
\end{Theorem}

\proof
Consider the following morphisms:
\begin{eqnarray}
\label{uni}
f^+: \overline{\mathcal{M}}_{0,6} &\longrightarrow& \overline{S}_{2}^+, \nonumber \\  & & \\
f^-: \overline{\mathcal{M}}_{0,6} &\longrightarrow&  \overline{S}_{2}^-. \nonumber
\end{eqnarray}

In order to define $f^+$ and $f^-$, let $(C; p_1, \ldots p_6)$ be a 
$6$-pointed, stable, genus zero curve. The morphism $f^+$ (respectively 
$f^-$) associates to $(C; p_1, \ldots p_6)$ the admissible covering $Y$ of $C$ which is branched at the 
$p_i$'s. Generically, it comes equipped with the line bundle $\mathcal{O}_Y(q_1+q_2-q_3)$  (respectively $\mathcal{O}_Y(q_1)$), where $q_i$ denotes the point of $Y$ lying 
above $p_i$, and the morphism extends to the boundary since 
$\overline{\mathcal{M}}_{0,6}$ is normal and $\overline{S}_{2}$ is finite 
over $\overline{\mathcal{M}}_{2}$.

By \eqref{uni}, we have
$$
H^k(\stwobar) \hookrightarrow H^k(\overline{\mathcal{M}}_{0,6}) .  
$$

The claim follows from \cite{Keel:92} since $H^k(\overline{\mathcal{M}}_{0,n})=0$ for every 
$n$ and every odd $k$.

\endproof

In order to determine $H^2(\stwobar)$, we first look at algebraic 
cycles of codimension one in $\stwobar$.
We recall that the boundary $\partial \sgbar = \sgbar \setminus S_g$ is 
the union of the irreducible components $A_i^+$, $B_i^+$ 
(contained in $\sgbar^+$) and $A_i^-$, $B_i^-$ (contained in $\sgbar^-$), 
which are completely described in \cite{Cornalba:89}, \S~7, by their 
general members:  
\begin{itemize}
\item $A_i^+$ (resp. $B_i^+$ ), $i > 0$: two smooth components $C_1$ and 
$C_2$ of genera $i$ and $g - i$, joined at points $p\in C_1$ and $q\in C_2$ 
by a $\P^1$, with even (resp. odd) theta-characteristics on $C_1$ and $C_2$ 
\item $A_i^-$  (resp. $B_i^-$), $i>0$: as above, but with an even (resp. odd) 
theta-characteristic on $C_1$ and an odd (resp. even) one on $C_2$
\item $A_0^+$ (resp. $A_0^-$): an irreducible curve of genus $g$ with only one 
node, with an even (resp. odd) spin structure
\item $B_0^+$ (resp. $B_0^-$): an irreducible curve of genus $g$ with
only one node, blown up at the node (so to add an exceptional
component $E \cong \P^1$), with an even (resp. odd) spin structure 
glued to ${\mathcal O}_E(1)$ on $E$.
\end{itemize}
Let $\alpha_i^+$, $\beta_i^+$, $\alpha_i^-$, $\beta_i^-$ denote the
corresponding classes in $A_{3g-4}(\sgbar)$ for $i=0,1$.  
In \cite{Cornalba:89}, p. 585, the class $\beta_1^-$ is defined to be zero, 
hence we always omit it.  We point out that the arguments provided in
\cite{Cornalba:89} in order to prove Proposition~7.2 imply that
$\alpha_i^+$, $\beta_i^+$, $\alpha_i^-$, $\beta_i^-$ are independent
for $g=2$, too.

\begin{Proposition}\label{bnd}
The Chow group $A_2(\stwobar)$ is generated by boundary classes.
\end{Proposition}

\proof By restricting $f^+$ (resp. $f^-$) over $S^+_2$ (resp. $S^-_2$),  we obtain finite surjective morphisms from $\mathcal{M}_{0,6}$ to $S_2$. By \cite{Keel:92},  $A_2(\overline{\mathcal{M}}_{0,6})$
is generated by boundary classes. Thus,  the exact sequence
$$
A_2(\overline{\mathcal{M}}_{0,6} \setminus \mathcal{M}_{0,6}) 
\to A_2(\overline{\mathcal{M}}_{0,6}) \to A_2(\mathcal{M}_{0,6})
$$
yields $A_2(\mathcal{M}_{0,6})=0$. By \cite{Faber}, 
Lemma A, we conclude that $A_2(S_2)=0$, so the map 
$A_2(\stwobar \setminus S_2) \to A_2(\stwobar)$ is onto.

\qed

Next, we are going to prove that the whole cohomology of $\stwobar$ is indeed algebraic. 
In order to do so, we need the following auxiliary result. 

\begin{Lemma}
\label{unilem}
$\overline{S}_{1,n}^{\hspace{0.05cm}(m_1, \ldots, m_n)}$
is a unirational variety for $n \le 2$. 
\end{Lemma}

\proof Since the structure sheaf is the unique odd theta-characteristic 
on an elliptic curve, we have $\overline{S}_{1,n}^{\hspace{0.05cm}(0, \ldots, 0), -} \cong \overline{\mathcal{M}}_{1,n}$ for every $n$, and 
${\overline{\mathcal M}}_{1,n}$ is rational for $n\leq 10$, 
as shown in \cite{Bel}.
Next, we are going to define three surjective morphisms by defining them generically 
as in the proof of Theorem~\ref{odd}. The first one is
$$
h_4: \overline{\mathcal{M}}_{0,4} \longrightarrow 
\overline{S}_{1,1}^{\hspace{0.05cm}(0), +}.
$$

Let $(C; p_1, p_2, p_3, p_4)$ be a $4$-pointed stable genus zero 
curve. 
If $E$ is the admissible covering of $C$ that is branched at the $p_i$'s and $q_i$ 
denotes the preimage of $p_i$ under such covering, then the morphism $h$ maps the pointed curve $(C; p_1, p_2, p_3, p_4)$  to $(E; q_1)$ with the even theta-characteristic $\mathcal{O}_E(q_1 - q_2)$. 

In a similar fashion, we define a second morphism as follows:

$$h_5: \overline{\mathcal{M}}_{0,5} \longrightarrow 
\overline{S}_{1,2}^{\hspace{0.05cm}(0,0), +}.
$$

Let  $(C; p_1, p_2, p_3,$ $p_4, p)$ be a $5$-pointed stable genus zero curve. The image under $h_5$ is the pointed curve $(F; r_1, r)$ with the even theta-characteristic  $\mathcal{O}_F(r_1 - r_2)$.  The curve $F$ is the admissible covering of $C$ that is branched at the points $p_i$'s. The point $r_i$ is the preimage of $p_i$, while the point $r$ is one of the two points in the preimage of $p$. Notice that a different choice will yield the same point in $\overline{S}_{1,2}^{\hspace{0.05cm}(0,0), +}$. Finally, we have a third surjective morphism 
$$
t_5: \overline{\mathcal{M}}_{0,5} \longrightarrow 
\overline{S}_{1,2}^{\hspace{0.05cm}(1,1)}, 
$$
which is defined as follows. With the same notation used for $h_5$, the image under 
$t_5$ of $(C; p_1, p_2, p_3,$ $p_4, p)$ is given by $(E; r_2, r_3)$ with the odd theta-characteristic $\mathcal{O}_E(r_1)$.  The points $r_2$ and $r_3$ lie in the preimage of $p$ under $t_5$. Note that the pointed curve $(E; r_3, r_2)$ and $\mathcal{O}_E(r_1)$ will yield the same element in $\overline{S}_{1,2}^{\hspace{0.05cm}(1,1)}$.

Since the moduli spaces $\mobar$ are rational, the claim is completely proved.

\endproof

\begin{Proposition}\label{hodge}
The Hodge diamond of $\stwobar$ is given by 
\smallskip
$$
\xymatrix{
 &  & & 2 \ar@{-}[dr] \ar@{-}[dl] & & &\\
& &  0 \ar@{-}[dr] \ar@{-}[dl] & & 0 \ar@{-}[dr] \ar@{-}[dl] & & \\
& 0 \ar@{-}[dr] \ar@{-}[dl] & & 7 \ar@{-}[dr] \ar@{-}[dl] & & 0 \ar@{-}[dr] \ar@{-}[dl] &\\
 0  & & 0  & & 0  & & 0  \\
& 0 \ar@{-}[ur] \ar@{-}[ul]& & 7 \ar@{-}[ur] \ar@{-}[ul] & & 0 \ar@{-}[ul] \ar@{-}[ur]\\
& & 0 \ar@{-}[dr]  \ar@{-}[ur] \ar@{-}[ul]& & 0 \ar@{-}[dl] \ar@{-}[ur]\ar@{-}[ul]& &\\
& & & 2 \ar@{-}[ur] \ar@{-}[ul] & & & 
}
$$
\centerline{{\sc Figure 3}: {\rm The Hodge Diamond of} $\stwobar$}
\end{Proposition}

\proof By Theorem \ref{odd} and Proposition \ref{bnd}, all we need to
show is that $h^{2,0}(\stwobar)=0$.  As recalled in Remark
\ref{remaff}, each connected component of $S_2$ is an affine variety. Hence,
$H_k(S_2)=0$ for $k \ge \dim(S_2)+1=4$. By Poincar\'{e} duality, we have
$H^k_c(S_2)=0$ for $0 \le k \le 2$. By the exact sequence
$$
\ldots \to H^k_c(S_2) \to H^k(\stwobar) \to H^k(\partial \stwobar) \to \ldots 
$$ 
we get an injective morphism $
H^k(\stwobar) \rightarrow H^k(\partial \stwobar)
$, which is compatible with the Hodge structures (exactly as 
in \cite{ArbCor:98}, p. 102); so there is an injection
$
H^{p,0}(\stwobar) \hookrightarrow H^{p,0}(\partial \stwobar) 
$
for $p \le 2$. Notice that each irreducible component of $\partial \stwobar$
is the image of a map from either ${\overline{S}}_{1,1} \times 
{\overline{S}}_{1,1}$ or ${\overline {S}}_{1,2}^{\hspace{0.05cm}(0, 0)}$ 
or ${\overline {S}}_{1,2}^{\hspace{0.05cm}(1, 1)}$. Analogously to Lemma 2.6 in \cite{ArbCor:98},  it suffices to check that 
$h^{p,0}(\overline{S}_{1,n}^{\hspace{0.05cm}(m_1, \ldots, m_n)})=0$ for $n \le 2$.  By Lemma \ref{unilem}, we get a dominant rational map
$
\P^{n} \dashrightarrow 
\overline{S}_{1,n}^{\hspace{0.05cm}(m_1, \ldots, m_n)}
$. Thus, as in \cite{GH}, p. 494, we have an injective morphism
$
H^{p,0}(\overline{S}_{1,n}^{\hspace{0.05cm}(m_1, \ldots, m_n)}) \hookrightarrow H^{p,0}(\P^{n})
$. In particular, $h^{2,0}(\stwobar)=0$, and the claim follows.

\qed

\begin{remark}\label{referee}
As pointed out by the referee, the vanishing of $h^{2,0}(\stwobar)$ is also 
a direct consequence of (\ref{uni}) and the rationality of 
$\overline{\mathcal{M}}_{0,6}$.
\end{remark}

\begin{Theorem} \label{h2}
The rational cohomology of $\stwobar$ is algebraic. In particular, the
following hold: \\ i) a basis for the cohomology group $H^2(\stwobar)$
is given by the boundary classes $\alpha_i^+$, $\alpha_i^-$,
$\beta_i^+$, $\beta_i^-$; \\ ii) a basis for the cohomology group
$H^4(\stwobar)$ is given by the products $\Delta \alpha_i^+$, $\Delta
\alpha_i^-$, $\Delta \beta_i^+$, $\Delta \beta_i^-$, where $\Delta$ is
the sum of all boundary classes.
\end{Theorem}

\proof By the standard exponential sequence, we have a long exact sequence
$$ \ldots \rightarrow H^1(\stwobar, {\mathcal O}) \otimes {\mathbb Q}
\rightarrow H^1(\stwobar, {\mathcal O}^*) \otimes {\mathbb Q} \rightarrow H^2(\stwobar,
{\mathbb Q}) \rightarrow H^2(\stwobar, {\mathcal O}) \otimes {\mathbb
Q} \rightarrow \ldots
$$

By Proposition \ref{hodge} we have $H^1(\stwobar, {\mathcal O})
\otimes {\mathbb Q} = (0) =H^2(\stwobar, {\mathcal O}) \otimes
{\mathbb Q}$. Hence the (rational) cohomology group $H^2(\stwobar)$
is isomorphic to $Pic (\stwobar) \otimes {\mathbb Q}$.  Thus, i)
follows from Proposition \ref{iso} and Proposition \ref{bnd}. As for
ii), note that $\Delta$ is ample because $S_2$ is affine
(cfr. Remark \ref{remaff}).  Accordingly, multiplication by
$\Delta$ induces an isomorphism between $H^2(\stwobar)$ and
$H^4(\stwobar)$ by the Hard Lefschetz Theorem, which holds in the
category of projective orbifolds - see \cite{Steenbrink:77},
Theorem~1.13. By Theorem \ref{odd}, the claim is completely proved.
\qed

\begin{Corollary}
$A^*(S^+_2) \cong A^*(S^-_2)\cong {\mathbb Q}$.
\end{Corollary}

\begin{remark}\label{genevodd}
Theorem \ref{h2} gives generators for the rational cohomology of $\stwobar^+$ (the classes with superscript $+$) and $\stwobar^-$ (the classes with superscript $-$).
\end{remark}

\begin{Corollary}
The Euler characteristic of  the normal variety $\stwobar$ is $18$. In particular, $e(\stwobar^+)=10$ and $e(\stwobar^-)=8$.
\end{Corollary}

For the ring structure of $H^*({\overline S}_2)$, we need some results on the cohomology of  genus one spin curves. Denote by $\alpha_0^+$, $\alpha_1^+$ (resp. $\alpha_0^-$,
$\alpha_1^-$) the (Poincar\'{e}) dual cohomology class of the loci $A_0^+$, $A_1^+$
(resp., $A_0^-$, $A_1^-$) in ${\overline S}_{1,n}$ described above in the case $n=0$
(the only difference here is that the rational component of the general member carries 
$n$ marked points). Next, define by $\lambda_1^+$ (resp. $\lambda_1^-$) the pull-back on 
$\overline{S}_{1,n}^+$ (resp. $\overline{S}_{1,n}^-$) of the cohomology class $\lambda$ on ${\overline {\mathcal M}}_{1,n}$. Clearly, these classes can be expressed in terms of
boundary divisors. In fact, we have $\lambda = \frac{1}{12} \delta_{irr}$ 
on ${\overline {\mathcal M}}_{1,n}$. Moreover, the space ${\overline S}^+_{1,1}$ is isomorphic to ${\mathbb P}^1$ and if $\gamma_{1,2}^+: {\overline S}^+_{1,2} \rightarrow {\overline {\mathcal M}}_{1,2}$ is the $3:1$ covering of ${\overline {\mathcal M}}_{1,2}$, then we have $(\gamma_{1,2}^+)^*(\delta_{irr})= 3\alpha_0^+$. Hence, 
\begin{equation}
\lambda_1^+ = \frac 14 \alpha_0^+, \qquad \lambda_1^- = \frac{1}{12} \alpha_0^-.
\label{elleuno}
\end{equation}

\begin{remark}\label{rem}
We have $(\lambda_1^+)^2=(\lambda_1^-)^2=0$. Indeed, as pointed out in \cite{Bel}, 
Lemma~3.2.6, $\delta_{irr}^2=0$ on ${\overline {\mathcal M}}_{1,n}$ since the 
boundary divisor $\Delta_{irr}$ is the pull-back of a point in 
${\overline {\mathcal M}}_{1,1} \cong \mathbb{P}^1$ via the natural 
forgetful morphism.  
\end{remark}

\begin{Theorem}
\label{coho2n}
The cohomology ring $H^*({\overline S}_2^+)$ is isomorphic to the quotient ring ${\mathbb Q}[\alpha_0^+, \, \alpha_1^+, \, \beta_0^+, \, \beta_1^+]/ J^+ $, where $J^+$ is the ideal generated by the following elements:  

$$\alpha_1^+\beta_1^+, \qquad \beta_0^+\beta_1^+, \qquad \alpha_0^+\alpha_1^+ - \beta_0^+\alpha_1^+, 
$$

$$ 
\alpha_0^+ \alpha_1^+  + 8\left( \alpha_1^+\right)^2, \qquad \alpha_0^+ \beta_1^+  \,+ 24 (\beta_1^+)^2, \qquad  4(\beta^+_0)^2+ 8\alpha_1^+ \beta_0^+ - 3\alpha_0^+ \beta_0^+ ,
$$ 

$$
(\alpha_0^+)^2\alpha_1^+, \qquad (\alpha_0^+ )^2\beta_0^+, \qquad 3(\alpha_0^+)^3 + 22(\alpha_0^+)^2\beta_1^+.
$$

\end{Theorem}

\proof The relations $\alpha_1^+\beta_1^+, \; \beta_0^+\beta_1^+$
follow immediately from the definition of $\beta_1^+$. In fact, the
intersection of the corresponding cycles are empty by parity reasons
(for instance, if $B_0^+ \cap B_1^+$ were non-empty, then it would be 
one-dimensional and its general element should carry an 
even theta characteristic restricting to an even one on the genus $0$ 
component and to an odd one to the genus $1$ component, which is 
clearly impossible). 

Let $\theta$ be the boundary map

\begin{equation}
\label{tetta}
\theta: {\overline {S}}^+_{1,1} \times {\overline {S}}^+_{1,1} \rightarrow {\overline{S}}_2^+
\end{equation}
defined as follows. A point in the domain is the datum of two  $1$-pointed genus one curves  with even spin structures. The image point is obtained by glueing the two curves along the two marked points and then blowing-up the node. The exceptional component $C$ is given a spin structure via ${\mathcal O}_C(1)$. We have
$$
\left(\alpha_0^+ - \beta_0^+\right)\alpha_1^+ = \frac 12 \theta_*\left(\theta^*\left(\alpha_0^+ - \beta_0^+\right)\right)=0
$$
since  $\theta^*\left(\alpha_0^+ - \beta_0^+\right)=0$. Indeed, the classes $\alpha_0^+$ and $\beta_0^+$ coincide on ${\overline {S}}^+_{1,1} \cong \mathbb{P}^1$, hence 
$1 \otimes \alpha_0^+ + \alpha_0^+ \otimes 1 - 1 \otimes \beta_0^+ - \beta_0^+ \otimes 1 =0$.

We recall from \cite{Mum}, p. 321, that $10\lambda= \delta_{irr} + 2\delta_1$ on ${\overline {\mathcal M}}_2$. We denote by $\lambda_2^+$ the pull-back of $\lambda$ on ${\overline{S}}_2^+$. By \cite{Cornalba:89}, Proposition~(7.2), i) and ii), we have 
$(\nu^+)^*(\delta_{irr})= \alpha_0^+  + 2\beta_0^+$ and 
$(\nu^+)^*(\delta_1)= 2\alpha_1^+ + 2\beta_1^+$, hence we obtain
\begin{equation}
\label{elle}
\lambda_2^+ = \frac{1}{10}(\alpha^+_0 + 2\beta_0^+ + 4\beta_1^+ + 4 \alpha_1^+).
\end{equation}

On the one hand, we get
$
\lambda_2^+ \alpha_1^+ = \frac{3}{10} \alpha_0^+\alpha^+_1 + \frac{2}{5}\left(\alpha_1^+\right)^2
$ 
by applying the previous relations.  On the other hand, if we recall (\ref{elleuno}) and apply the push-pull formula we see that 
\begin{eqnarray}
\lambda_2^+ \alpha_1^+ &=& \frac 12  \theta_*\left(\theta^*(\lambda_2)\right)= \frac 12 \theta_*\left( \lambda^+_1 \otimes 1  + 1 \otimes \lambda^+_1 \right) \\ &=& \frac 18 \alpha_0^+ \alpha_1^+ \nonumber.
\label{first}
\end{eqnarray}

Thus, the codimension two relation $\alpha_0^+ \alpha_1^+  + 8\left( \alpha_1^+\right)^2=0$ holds. 
 
As shown in \cite{Mum}, p. 321, the following relation holds in $H^4({\overline {\mathcal M}}_2)$:
 \begin{equation}
 \label{second}
 \lambda \delta_1 = \frac{1}{12} \delta_{irr} \delta_1.
 \end{equation}

By subtracting \eqref{first}, the pull-back of \eqref{second} yields
$
\lambda_2^+ \beta_1^+ = \frac{1}{12} \alpha_0^+\beta_1^+
$. By \eqref{elle}, the class $\lambda_2^+ \beta_1^+$ is also equal to  $1/10 \left( \alpha_0^+ \beta_1^+ + 4 \alpha_1^+ \beta_1^+ + 4 (\beta_1^+)^2 \right)$; so the codimension two relation 
\begin{equation}
\label{almost}
\alpha_0^+ \beta_1^+ + 24 (\beta_1^+)^2=0
\end{equation}
follows.  Notice that 
 $$
 3\left(\alpha_0^+ \alpha_1^+  + 8\left( \alpha_1^+\right)^2    \right) + \alpha_0^+ \beta_1^+ + 24 (\beta_1^+)^2=0
$$
is the pull-back of  $\delta_{irr} \delta_1 + 12 \delta_1^2$ on ${\overline {\mathcal M}}_2$, which is zero - see \cite{Mum}, p. 321.

In order to compute the last relation, let $\eta$ be the boundary map
$$
\eta: {\overline S}_{1,2}^+ \rightarrow {\overline S}^+_2
$$
which is defined as follows. A point in the domain is the datum of a two pointed genus one curve with an even spin structure. The corresponding image point is obtained by glueing the two marked points and blowing-up the node. The resulting quasi-stable curve is given a spin structure via the spin structure of the domain curve plus ${\mathcal O}_F(1)$ on the exceptional component $F$. By definition we have $\beta_0^+ = \frac 14 \eta_*(1)$. Notice that the order of the automorphism group of a generic curve in the image of $\eta$ is four. By \eqref{elleuno}, we have
\begin{equation}
\label{etarelation}
\lambda_2^+ \beta_0^+ = \frac 14 \eta_* \left ( \lambda_1^+\right) = \frac 14 \eta_* \left(\frac 14 \alpha_0^+\right) = \frac 14 \alpha_0^+ \beta_0^+ .
\end{equation}

By \eqref{elle}, the class $\lambda_2^+\beta^+_0$ is also equal to $1/10 \left( \alpha_0^+ \beta_0^+ + 2(\beta_0^+)^2  + 4 \alpha_1^+ \beta_0^+ \right)$. 
Hence, the last codimension two relation follows.

Let us now determine the codimension three relations. By \eqref{etarelation}, we get 
$$
(\lambda_2^+)^2\beta_0^+= \frac 14 \lambda_2^+ \beta_0^+ \alpha_0^+ = \frac{1}{16}(\alpha_0^+)^2(\beta_0^+).
$$

We remark that $(\lambda_2^+)^2 \alpha_0^+ = (\lambda_2^+)^2 \beta_0^+=0$. Indeed, restricting $(\lambda_2^+)^2$ to the divisors corresponding to $\alpha_0^+$ and $\beta_0^+$ is equivalent to computing $(\lambda_1^+)^2$ on ${\overline S}^+_{1,2}$, which is zero by Remark~\ref{rem}. This yields $(\alpha_0^+)^2(\beta_0^+)=0$. By plugging \eqref{elle} in this last identity, we get $(\alpha_0^+)^2\alpha_1^+=0$.  This, together with the codimension two relations above, yields that all degree three monomials in $H^6({\overline S}_2^+)$ vanish except for $
(\alpha_0^+)^3, \; (\alpha_0^+)^2\beta_1^+, \; \alpha_0^+(\beta_1^+)^2$. By \eqref{almost}, we have $(\beta_1^+)^2\alpha_0^+ = -\frac{1}{24}\alpha_0^+ \beta_1^+$. Finally, by the relation $(\lambda_2^+)^2\alpha_0^+=0$, we get $3(\alpha_0)^3 + 22(\alpha_0^+)^2\beta_1^+=0$. 

A calculation with computer program Macaulay \cite{BS} yields that the ring ${\mathbb Q}[\alpha_0^+, \alpha_1^+, \beta_0^+, \beta_1^+]/J^+$ is finite dimensional and has Betti numbers $1, 4, 4, 1$. Since this ring surjects onto $H^*({\overline S}^+_2)$ and the Betti numbers coincide, they are isomorphic. 

Notice in particular that $(\alpha_0^+)^2\beta_1^+= \frac{5}{4}[p]$, where $[p]$ is the class of a point on ${\overline S}^+_2$. This follows from the intersection number $ \delta_1^3= \frac{1}{576}$ on ${\overline {\mathcal M}}_2$.

\qed

\begin{Theorem}
\label{coho3n}
The cohomology ring  $H^*({\overline S}_2^-)$ is isomorphic to the quotient ring ${\mathbb Q}[\alpha_0^-, \alpha_1^-, \beta_0^-]/ J^- $, where $J^-$ is the ideal generated by the following  polynomials:

$$
3 (\beta^-_0)^2  + 6\alpha_1^- \beta_0^- - \alpha_0^- \beta_0^-,  \qquad  2 \alpha_1^- \beta_0^- - \alpha_1^ - \alpha_0^- , \qquad 12(\alpha_1^-)^2 + \alpha_1^-\alpha_0^-,
$$

$$
(\alpha_0^-)^2\beta_0^- \qquad 3 (\alpha_0^-)^3+ 32 \alpha_1^- (\alpha_0^-)^2. 
$$
\end{Theorem}

\proof  Analogously to (\ref{etarelation}), we can prove that 
\begin{eqnarray}
\label{orace}
\lambda_2^- \beta_0^- = \frac 16 \alpha_0^- \beta_0^- 
\end{eqnarray}
where $\lambda_2^-$ is the pull-back of $\lambda$ to ${\overline S}_2^-$, 
hence we have
$$
\lambda_2^-=\frac{1}{10}\left( \alpha_0^-+2\beta_0^-+4\alpha_1^-\right). 
$$

If we multiply by $\beta_0^-$ and compare with \eqref{orace}, we get the relation $3(\beta^-_0)^2  + 6\alpha_1^- \beta_0^- - \alpha_0^- \beta_0^-=0$ holds.  

Similarly to \eqref{tetta}, define the boundary map $\iota: {\overline {S}}^+_{1,1} \times {\overline {S}}^-_{1,1} \rightarrow {\overline{S}}_2^-$. 
Since $\alpha_0^+ = \beta_0^+$ on ${\overline {S}}^+_{1,1} \cong \mathbb{P}^1$,
we have:
\begin{eqnarray*}
\beta_0^-\alpha_1^-   &=&  \iota_* \iota^* \left( \beta_0^-\right)= \iota_*\left( \beta_0^+ \otimes 1 \right)  = \iota_*\left(\alpha_0^+ \otimes 1\right) \\ &=& \iota_*\left(\iota^*\left(\alpha_0^-\right)\right) - \iota_*\left( 1 \otimes \alpha_0^-\right)  = \alpha_0^- \alpha_1^-  - \iota_*\left(1 \otimes \alpha_0^-\right) \\ &=& \frac12 \alpha_0^- \alpha_1^-.
\end{eqnarray*}

Note that $\iota_*(1 \otimes \alpha_0^-) = \frac{1}{2} \alpha_0^- \alpha_1^-$. In fact,  the automorphism group of the generic curve in the intersection of the divisors corresponding to $\alpha_0^-$ and $\alpha_1^-$ has order two. Hence, the codimension two relation $2 \alpha_1^- \beta_0^- - \alpha_0^- \alpha_1^-=0 $ follows.  Finally, if we pull-back $\delta_{irr}\, \delta_1 + 12 \delta_1^2$ on ${\overline S}_2^-$, we get $48(\alpha_1^-)^2+ 2\alpha_1^- \alpha^-_0 + 4 \alpha_1^- \beta_0^-$, which yields the relation $12(\alpha_1^-)^2 + \alpha_1^-\alpha_0^-=0$.

If we use the codimension two relations, we obtain the following relations in codimension three:
$$
144 (\alpha_1^-)^3 - \alpha_1^- (\alpha_0^-)^2=0, \qquad 54(\beta_0^-)^3 - 6(\alpha_0^-)^2\beta_0^- + 45\alpha_1^-(\alpha_0^-)^2 =0, 
$$

$$
12 (\alpha_1^-)^2\alpha_0^-+ \alpha_1^- (\alpha_0^-)^2=0, \qquad 24(\alpha_1^-)^2\beta_0^- + \alpha_1^- (\alpha_0^-)^2=0,
$$

$$
\alpha_0^-(\beta_0^-)^2 + \alpha_1^- (\alpha_0^-)^2=0, \qquad 4(\beta_0^-)^2 \alpha_1^- - (\alpha_0^-)^2\alpha_1^-=0,
$$

$$ 2 \alpha_1^- \alpha_0^- \beta_0^- - \alpha_1^- (\alpha_0^-)^2=0.
$$

\medskip

By Remark~\ref{rem}, we have $(\lambda_2^-)^2\beta_0^-=(\lambda_2^-)^2\alpha_0^-=0$. By \eqref{orace}, we get $(\alpha_0^-)^2\beta_0^-=0$. If we plug in $ \lambda_2^- =  \frac{1}{10}\left(\alpha_0^- + 2 \beta_0^- + 4 \alpha_1^-\right)$ in $(\lambda_2^-)^2\alpha_0^-=0$ and use the relations above, we also have 
$
3 (\alpha_0^-)^3+ 32 \alpha_1^-(\alpha_0^-)^2=0.
$

A calculation with computer program Macaulay \cite{BS} yields that the quotient ring ${\mathbb Q}[\alpha_0^-, \alpha_1^-, \beta_0^-]/J^-$ is finite dimensional and has Betti numbers $1, 3, 3, 1$. Since this ring surjects onto $H^*({\overline S}^-_2)$ and the Betti numbers coincide, they are isomorphic. 

Finally, notice that   $\alpha_1^-(\alpha_0^-)^2= \frac{3}{16}[q]$, where $[q]$ is the class of a point on ${\overline S}_2^-$. This follows easily from the intersection number  $\delta_{irr} \delta_1^2= -\frac{1}{48}$ on ${\overline {\mathcal M}}_2$. 
\qed

\centerline{\sc appendix A}

\medskip

We recall that the moduli space of stable curves has a natural
stratification by topological type. An explicit description in genus $2$ 
can be found in \cite{Mum}: there are seven strata $\Delta_{G_i}$,
which correspond to the graphs $G_i$ in Figure 1.

\begin{figure}[h]
\begin{center} 
\mbox{\epsfig{file=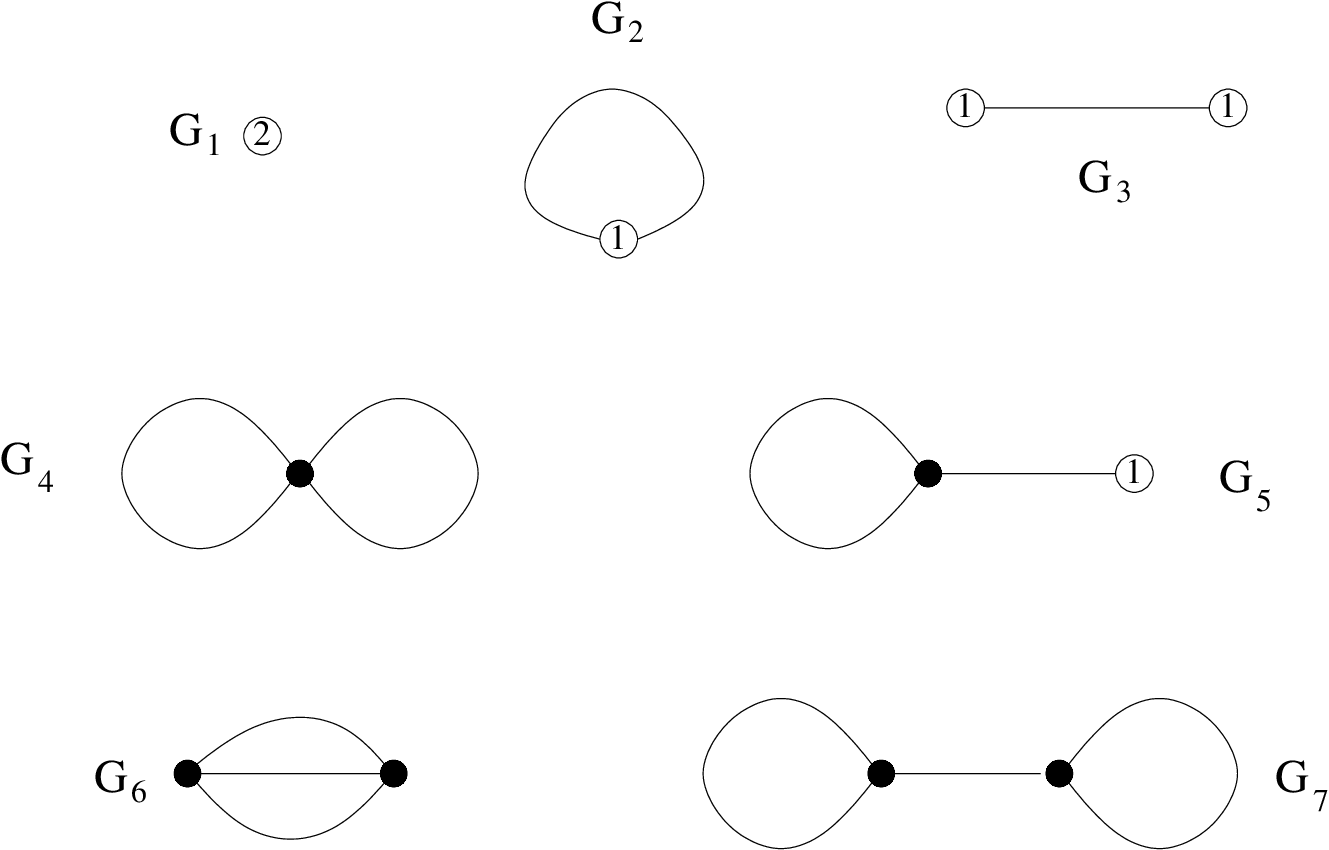,width=9.5cm,height=6cm,angle=0}}
\caption[]{Graph-types.}
\label{types}
\end{center}
\end{figure}

\bigskip

We shall list the components of $\nu^{-1}\left(\Delta_{G_i}\right)$ by describing their general member.

The preimage of $\Delta_{G_1}$ yields two strata, namely $S_g^+$ and $S^-_g$. 

\smallskip

The eight strata of the preimages $\nu^{-1}(\Delta_{G_2})$ and  $\nu^{-1}(\Delta_{G_3})$ are described in \cite{Cornalba:89}.  Their closures are denoted by $A^+_j, A^-_j, B^+_j, B^-_j$.

\smallskip

The preimage $\nu^{-1}(\Delta_{G_4})$ is the union of the following strata: \\
$\bullet$ $C^+$: an irreducible curve with two nodes and an even spin structure; \\
$\bullet$ $C^-$: an irreducible curve with two nodes and an odd spin structure; 

\begin{figure}[h]
\begin{center} 
\mbox{\epsfig{file=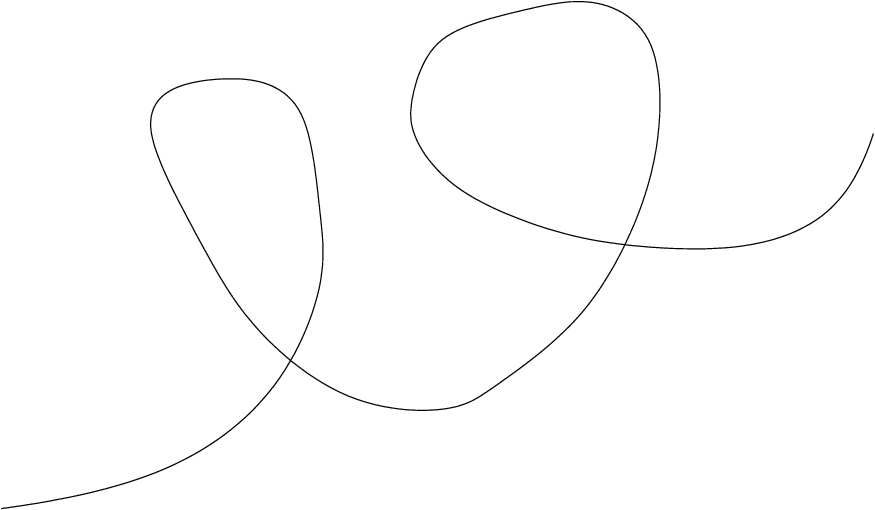,width=4cm,height=2.5cm,angle=0}}
\caption[]{Curves in $C^+$ and $C^-$.}
\end{center}
\end{figure}
\noindent $\bullet$ $D^+$: a rational irreducible curve with two nodes - blown-up at one of them - with an even spin structure; \\
$\bullet$ $D^-$: as above, but with an odd spin structure; 

\begin{figure}[h]
\begin{center} 
\mbox{\epsfig{file=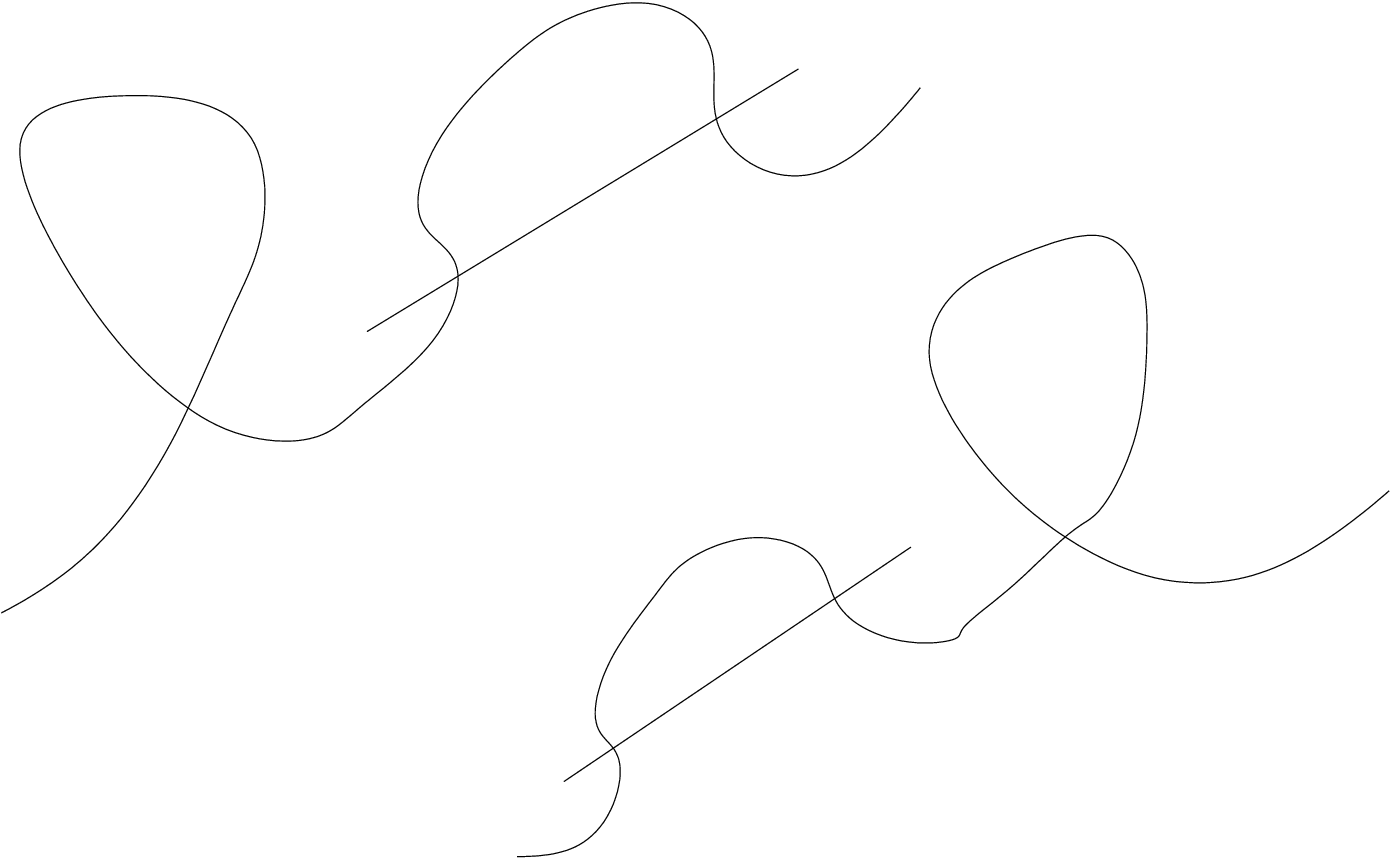,width=5cm,height=3cm,angle=0}}
\caption[]{Curves in $D^+$ and $D^-$.}
\end{center}
\end{figure}
\noindent $\bullet$ $E$: a rational irreducible curve with two nodes - blown-up at both nodes - with an even spin structure. 

\begin{figure}[h]
\begin{center} 
\mbox{\epsfig{file=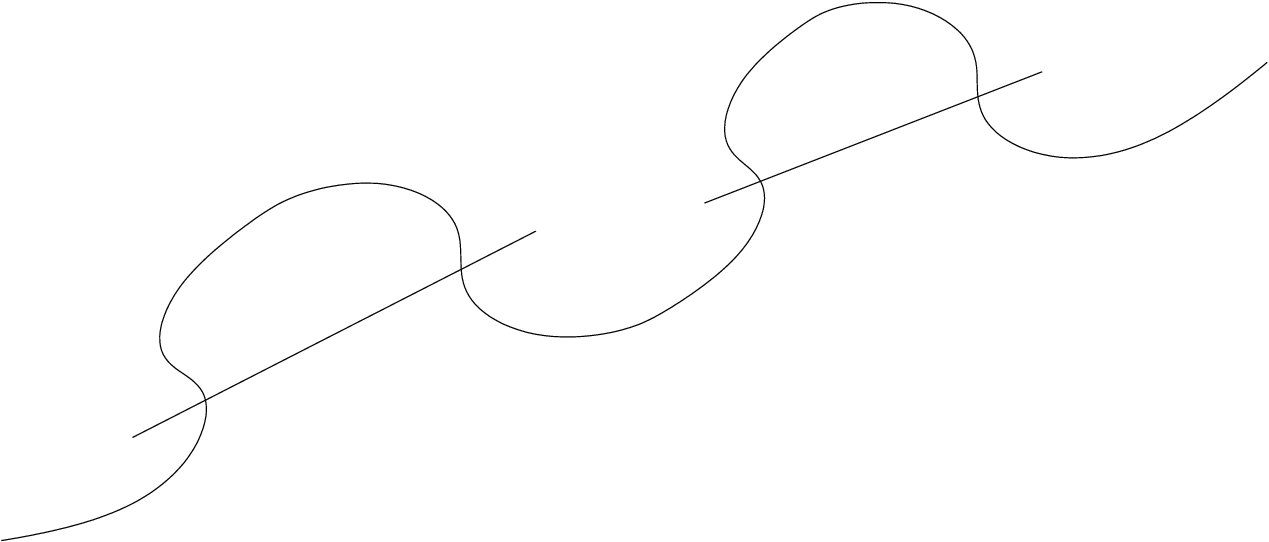,width=5cm,height=3cm,angle=0}}
\caption[]{Curves in $E$.}
\end{center}
\end{figure}

The preimage $\nu^{-1}(\Delta_{G_5})$ is the disjoint union of the following strata: \\
$\bullet$ $X^+$: a smooth elliptic curve $F_1$ and an irreducible nodal curve $F_2$ - joined via a smooth rational curve at two points $y_1$ and $y_2$ - with even theta-characteristics on $F_1$ and $F_2$; \\
$\bullet$ $X^-$: as above, but with an even theta-characteristic on $F_1$ and an odd theta-characteristic on $F_2$; \\
$\bullet$ $Y^+$: as above, but with odd theta-characteristics on $F_1$ and $F_2$; \\
$\bullet$ $Y^-$: as above, but with an odd theta-characteristic on $F_1$ and an even theta-characteristic on $F_2$; 

\begin{figure}[h]
\begin{center} 
\mbox{\epsfig{file=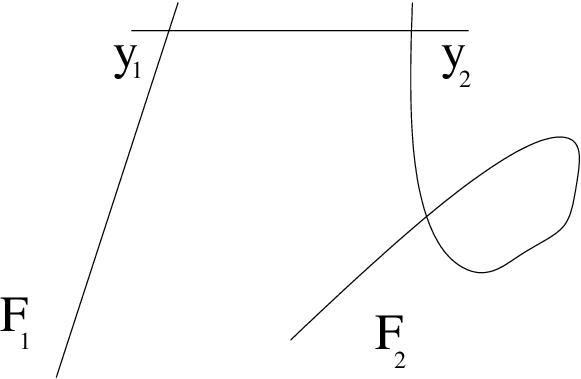,width=5cm,height=3cm,angle=0}}
\caption[]{Curves in $X^+$, $X^-$, $Y^+$ and $Y^-$.}
\end{center}
\end{figure}
\noindent $\bullet$ $Z^+$: an irreducible nodal curve $T_2$ (blown-up at the node) and a smooth elliptic curve $T_1$ - joined via a smooth rational curve - with an even theta-characteristic on $T_1$ and an even spin structure on $T_2$; \\
$\bullet$ $Z^-$: as above, but with an odd theta-characteristic on $T_1$ and an even spin structure on $T_2$.

\begin{figure}[h]
\begin{center} 
\mbox{\epsfig{file=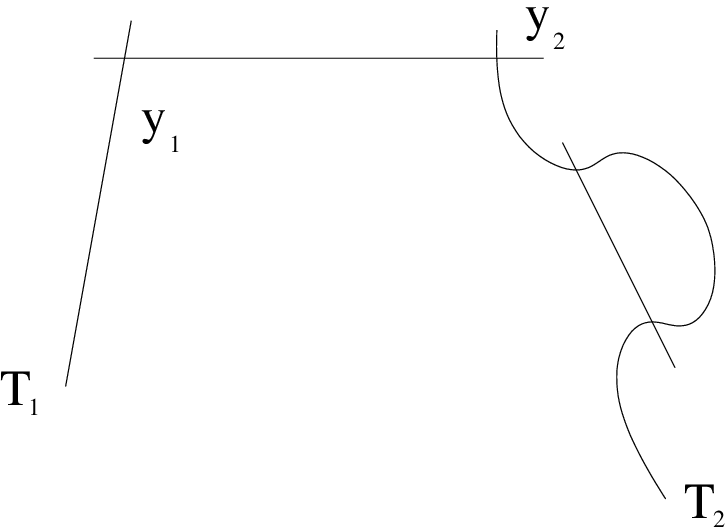,width=5cm,height=3cm,angle=0}}
\caption[]{Curves in $Z^+$ and $Z^-$.}
\end{center}
\end{figure}

\smallskip

The preimage $\nu^{-1}(\Delta_{G_6})$ has seven strata, i.e., \\
$\bullet$ $L^+$: two smooth rational curves that intersect in three points (blown-up at one of them) with an even spin structure\footnote{Notice that blowing-up at a different point yields an isomorphic spin curve.}; \\
$\bullet$ $L^-$: as above, but  with an odd spin structure; 

\begin{figure}[h]
\begin{center} 
\mbox{\epsfig{file=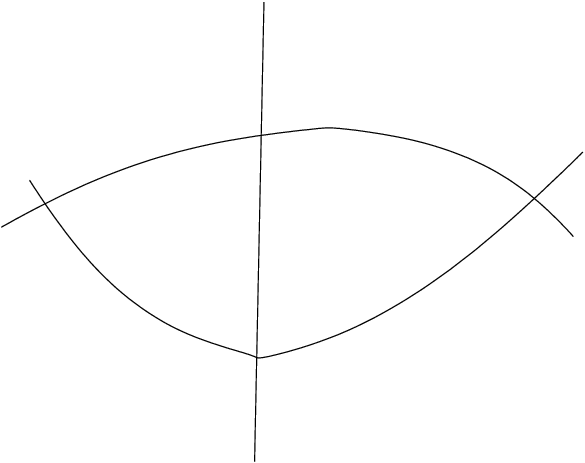,width=5cm,height=3cm,angle=0}}
\caption[]{Curves in $L^+$ and $L^-$.}
\end{center}
\end{figure}
\noindent $\bullet$ $M$: two smooth rational curves $K_1$ and $K_2$ that intersect in three points $p_1, p_2, p_3$ (blown-up at the intersection points) with an even spin structure.

\begin{figure}[h]
\begin{center} 
\mbox{\epsfig{file=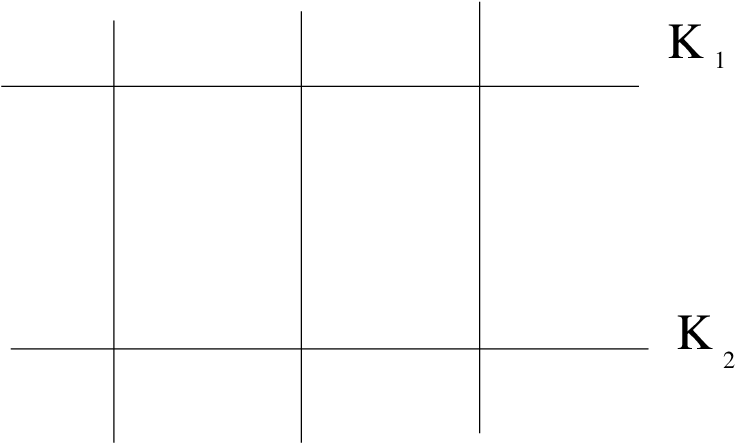,width=5cm,height=3cm,angle=0}}
\caption[]{Curves in $M$.}
\end{center}
\end{figure}

\smallskip

The preimage $\nu^{-1}(\Delta_{G_7})$ has seven strata (points), i.e., \\
$\bullet$ $P^+$: two irreducible nodal curves $C_1$ and $C_2$ of (arithmetic) genus $1$ - joined at points $p$ and $q$ by a smooth rational curve - with an even theta-characteristics on $C_1$ and $C_2$; \\
$\bullet$ $Q^+$: as above, but with odd theta-characteristics on $C_1$ and $C_2$; \\
$\bullet$ $P^-=Q^-$: as above, but with an odd theta-characteristic on $C_1$ and an even theta-characteristic on $C_2$;\\
\begin{figure}[h]
\begin{center} 
\mbox{\epsfig{file=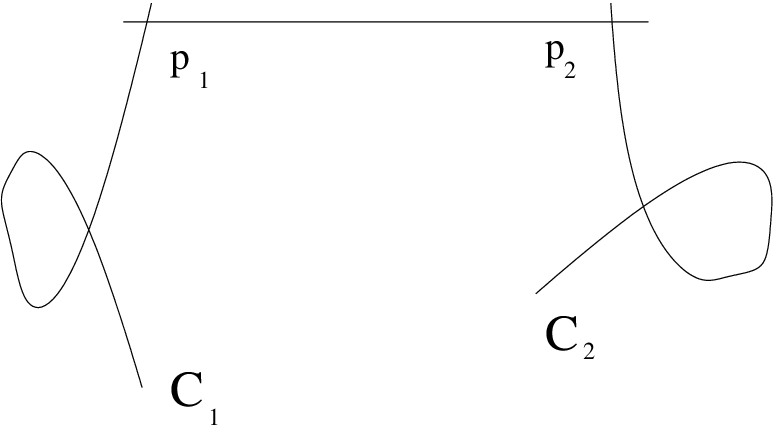,width=5cm,height=3cm,angle=0}}
\caption[]{Curves in $P^+, Q^+$ and $P^-=Q^-$.}
\end{center}
\end{figure}
\noindent $\bullet$ $R$: two irreducible curves with one node (each of them blown-up at the node) that are joined by a smooth rational curve at two distinct points, with an even spin structure; 

\begin{figure}[h]
\begin{center} 
\mbox{\epsfig{file=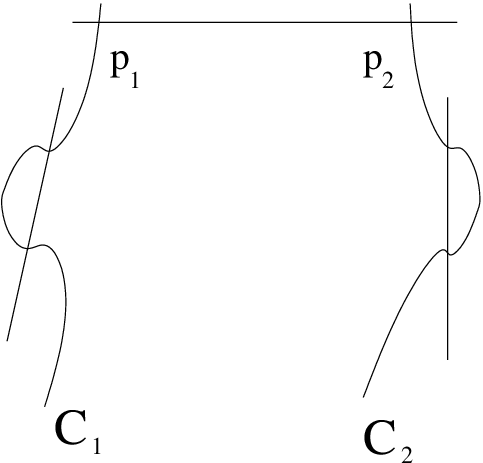,width=5cm,height=3cm,angle=0}}
\caption[]{Curves in $R$.}
\end{center}
\end{figure}
\noindent $\bullet$ $U^+$: an irreducible curve $V_2$ with one node (blown-up at the node) and an irreducible nodal curve $V_1$ of arithmetic genus one - joined  via a smooth rational curve at two distinct points -  with an even theta-characteristic on $V_1$ and an even spin structure on $V_2$;\\
$\bullet$ $U^-$: as above, but with an odd theta-characteristic on $V_1$ and an even spin one on $V_2$.

\begin{figure}[h]
\begin{center} 
\mbox{\epsfig{file=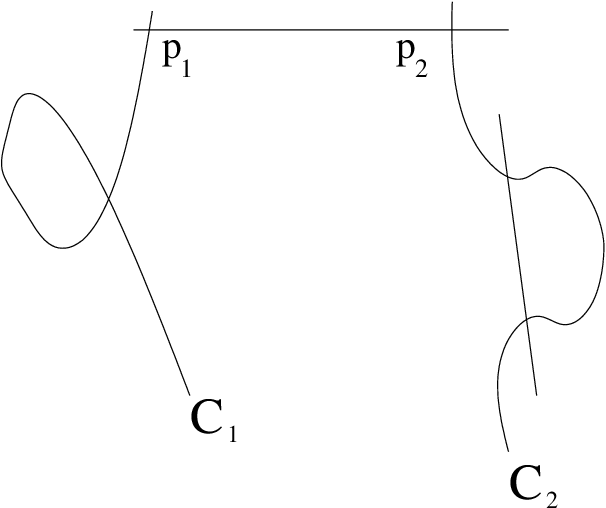,width=5cm,height=3cm,angle=0}}
\caption[]{Curves in $U^+$ and $U^-$.}
\end{center}
\end{figure}

\end{document}